\documentclass{article}
\usepackage[T2A]{fontenc}
\usepackage{amssymb,amsmath,amsthm}
\usepackage{hyperref}

\newcommand{\Q}{\mathbb{Q}}
\newcommand{\Z}{\mathbb{Z}}

\title{On strongly regular graph with parameters $(65,\,32,\,15,\,16)$}
\author{Oleg Gritsenko}
\begin{document}

\maketitle
  \begin{abstract}
    We construct a strongly regular graph with the parameters $(65,\,32,\,15,\,16)$. The idea is to
    search for an adjacency matrix that consists of circulant blocks. Equations with such matrices
    can be reduced to congruences with polynomials matrices of smaller orders. We can consider these
    congruences over different moduli for a more efficient computational approach.
  \end{abstract}
  
  \section{Introduction}
  
  A strongly regular graph (SRG) with parameters $(v,\,k,\,\lambda,\,\mu)$ is an undirected graph with
  $v$ vertices, each one has $k$ edges and every two adjacent or non-adjacent vertices have $\lambda$
  or $\mu$ common neighbours respectively. A special subset of SRG are
  conference graphs with parameters $k = \frac{v-1}2$, $\lambda = \frac{v-5}4$ and
  $\mu = \frac{v-1}4$. Conference graphs are related to symmetric conference matrices, i.e., matrices
  $C$ with 0 on the diagonal and $\pm1$ off the diagonal such that $C^2 = (n-1)I_n$, where $n$ is the
  order of $C$ and $I_n$ is the identity matrix of order $n$.
  
  An existence of a SRG for given parameters is a question of a great interest.
  There are several necessary conditions, however, for many parameters it is unknown whether there
  is such a SRG. All feasible sets with $v\leqslant 1300$ with constructions for known graphs are
  tracked at Andries Brouwer's website \cite{brouwer}. Until recently, the smallest unknown case was
  a conference graph with the parameters $(65,\,32,\,15,\,16)$.
  
  Denote $J_n$ a matrix of order $n$ consisting of ones, $c_n(x) = 1+x+x^2+\ldots+x^{n-1}$.
  It is well known that if $A$ is an adjacency matrix of a SRG with the parameters
  $(v,\,k,\,\lambda,\,\mu)$, it must satisfy the equations
  
  \begin{equation}
    \label{main_eq}
    A^2 + (\mu - \lambda) A = (k - \mu)I_n + \mu J_n,\;\;\;\; AJ_n = k J_n,
  \end{equation}
  where $n = v$. Also denote $M_n(R)$ a ring of matrices of order $n$ over a ring $R$.
  
  \section{Circulant matrices and matrix blocks}

  Suppose we a trying to find some solution of a matrix equation, but iterating over all possible
  matrices would take too much time. We can try to search for a matrix with some specific
  structure, for example, we can suppose it consists of circulant blocks. In this section we'll see how
  to represent such matrices as polynomial matrices of a smaller order. It might be easier
  to search for solutions having this specific form, but there is no guarantee that such a solution
  would exist, even if there is a solution of the original equation.
  
  Suppose we have two circulant matrices $B,\,C \in M_m(\Z)$, i.e. $B_{ij} =
  b_{(i-j)\!\!\! \mod m}$ and $C_{ij} = c_{(i-j)\!\!\! \mod m}$. There is a one-to-one
  correspondence between such matrices and polynomials in $\Z[x]/(x^m-1)$, namely, $B$ and $C$
  correspond to $b(x) = \sum\limits_{i=0}^{m-1} b_i x^i$ and $c(x) = \sum\limits_{i=0}^{m-1} c_i x^i$
  respectively. Note that the product $BC$ is also a circulant matrix with the corresponding
  polynomial $b(x) c(x) \pmod{x^m-1}$.

  Now, if a matrix $B \in M_{lm}(\Z)$ can be split into $l\times l$ blocks, where each block
  $B_{ij} \in M_m(\Z)$ is circulant, we consider the matrix $B(x) \in M_l(\Z[x] / (x^m-1))$
  consisting of the polynomials $b_{ij}(x)$ corresponding to these blocks $B_{ij}$. If another
  matrix $C \in M_{lm}(\Z)$ has the same structure and corresponds to $C(x) \in M_l(\Z[x] / (x^m-1))$,
  their product $BC$ also has the same structure, i.e., consists of $l\times l$ circulant blocks,
  and its corresponding polynomial matrix is $B(x)\cdot C(x) \pmod{x^m-1}$.
  
  For example, consider the adjacency matrix $B$ of the Petersen graph, a SRG with the parameters
  $(10,\,3,\,0,\,1)$.

  $$
    \left(
    \begin{array}{ccccc|ccccc}
      0 & 1 & 0 & 0 & 1 & 1 & 0 & 0 & 0 & 0 \\
      1 & 0 & 1 & 0 & 0 & 0 & 1 & 0 & 0 & 0 \\
      0 & 1 & 0 & 1 & 0 & 0 & 0 & 1 & 0 & 0 \\
      0 & 0 & 1 & 0 & 1 & 0 & 0 & 0 & 1 & 0 \\
      1 & 0 & 0 & 1 & 0 & 0 & 0 & 0 & 0 & 1 \\
      \hline
      1 & 0 & 0 & 0 & 0 & 0 & 0 & 1 & 1 & 0 \\
      0 & 1 & 0 & 0 & 0 & 0 & 0 & 0 & 1 & 1 \\
      0 & 0 & 1 & 0 & 0 & 1 & 0 & 0 & 0 & 1 \\
      0 & 0 & 0 & 1 & 0 & 1 & 1 & 0 & 0 & 0 \\
      0 & 0 & 0 & 0 & 1 & 0 & 1 & 1 & 0 & 0 \\
    \end{array}
    \right)
  $$
  
  Note that $B^2 + B = 2I_{10} + J_{10}$, and that $B$ consists of 4 blocks $5\times5$, each one is a
  circulant matrix corresponding to the polynomials $x+x^4$, $1$, $1$, and $x^2+x^3$ respectively.
  So the matrix $B$ can be `compacted' to
  $$
    B(x) =
    \begin{pmatrix}
      x+x^4 & 1 \\
      1 & x^2 + x^3 \\
    \end{pmatrix}
  $$
  such that $B^2(x) + B(x) \equiv 2I_2 + c_5(x) J_2 \pmod{x^5-1}$.
  
  Note that if a `compacted' matrix can be split into circulant blocks too, we can repeat this
  procedure, using a different polynomial variable. To illustrate this, consider the
  Hoffman–Singleton graph, a SRG with the parameters $(50,\,7,\,0,\,1)$. For brevity we won't write
  its adjacency matrix here, but note it can be compacted to
  
  $$
    B(x,\,y) =
    \begin{pmatrix}
      x+x^4 & 1 + x (y + y^4) + x^4 (y^2 + y^3) \\
      1 + x^4 (y + y^4) + x (y^2 + y^3) & x^2 + x^3 \\
    \end{pmatrix}
  $$
  such that $B^2(x,\,y) + B(x,\,y) \equiv 6I_2 + c_5(x)c_5(y)J_2 \pmod{x^5-1,\,y^5-1}$. Taking its
  $y$-coefficients, we can expand it to the larger matrix $B(x) \in M_{10}(\Z[x] / (x^5-1))$
  such that $B^2(x) + B(x) \equiv 6I_{10} + c5(x)J_{10} \pmod{x^5-1}$. $B(x)$ in turn can be expanded
  to $B \in M_{50}(\Z)$ such that $B^2 + B = 6I_{50} + J_{50}$. Since $B$ is
  also a symmetric matrix consisting of zeros and ones and $BJ_{50} = 7J_{50}$, it is the adjacency
  matrix of a SRG with the parameters $(50,\,7,\,0,\,1)$.
  
  So if we a searching for some solution of some matrix equation, say, $A^2 + pA = qI_n + rJ_n$,
  where $n = ml$, we can try to search for a matrix consisting of $l\times l$ circulant blocks, and
  to reduce this matrix equation to $A^2(x) + pA(x) \equiv qI_l + r\,c_m(x)J_l \pmod{x^m-1}$. We
  can also consider this congruence modulo some factors of $x^m-1$, for example, $x-1$.
  
  Of course if we want $A$ to consist only of zeros and ones, polynomial coefficients of $A(x)$
  must also be zeros and ones. Furthermore, if we are looking for a symmetric matrix, $A(x)$ must
  satisfy the congruence $A^T(x) \equiv A(x^{m-1}) \pmod{x^m-1}$, in other words, $A_{ji}(x) \equiv
  A_{ij}(x^{m-1})\pmod{x^m-1}$.
  
  \section{Searching for $srg(65,\,32,\,15,\,16)$}
  
  Now let's search for a SRG with the parameters $(65,\,32,\,15,\,16)$. Its adjacency matrix
  $A \in M_{65}(\Z)$ must satisfy the equations $A^2 + A = 16I_{65} + 16J_{65}$ and
  $AJ_{65} = 32J_{65}$, according to (\ref{main_eq}).
  
  Keeping the previous section in mind, we can suppose that $A$ can be split into $5\times5$
  circulant blocks; or in $13\times13$ circulant blocks; or that the whole $A$ is circulant.
  Unfortunately, neither approach works.
  
  In the first case $A$ is reduced to a matrix $A_5(x) \in M_5(\Z[x] / (x^{13}-1))$ such that
  $$
    A^2_5(x) + A_5(x) \equiv 16I_5 + 16c_{13}(x) J_5 \pmod{x^{13}-1}.
  $$
  By taking $x = \zeta_{13} = e^{2\pi i/13}$, we get $A^2_5(\zeta_{13}) + A_5(\zeta_{13}) = 16I_5$. The
  eigenvalues of $A_5(\zeta_{13})$ are the solutions of the equation $\lambda^2 + \lambda = 16$,
  i.e., $\frac{-1\pm\sqrt{65}}2$. Thus $tr(A_5(\zeta_{13})) = k \frac{-1+\sqrt{65}}2 +
  (5-k) \frac{-1-\sqrt{65}}2$, where $0\leqslant k\leqslant 5$, and for any such $k$
  $tr(A_5(\zeta_{13})) \in \Q(\sqrt{65})\setminus \Q$. On the other hand, $tr(A_5(x)) \in \Z[x]$,
  so $tr(A_5(\zeta_{13})) \in \Z[\zeta_{13}]$. However, one can easily check that $\sqrt{65}
  \notin\Q(\zeta_{13})$.
  
  Similarly, $A$ can't consist of $13\times13$ circulant blocks, as $\sqrt{65}\notin\Q(\zeta_5)$.
  If $A$ itself is circulant, it corresponds to a polynomial $a(x) \in \Z[x]/(x^{65}-1)$, where
  $a^2(x) + a(x) \equiv 16 + 16c_{65}(x) \pmod{x^{65}-1}$. Once again we can take 
  $x = \zeta_5$ and conclude there is no such matrix.
  
  Instead, let's try the following modification. The first row of $A$ has 32 ones, so without loss
  of generality we may assume $A_{1j}=1$ for $2\leqslant j\leqslant33$ and $A_{1j}=0$ for
  $34\leqslant j\leqslant65$. Denote $B=(1\,1\,\ldots\,1\,0\,0\,\ldots\,0) = (A_{1j})$,
  $2\leqslant j\leqslant65$, $C=(A_{ij})\in M_{64}(\Z)$, $i,j\geqslant2$. Then the matrix $A$ can be
  written as
  $$
    \begin{pmatrix}
      0 & B \\
      B^T & C \\
    \end{pmatrix}
  $$
  The equality $A^2 + A = 16I_{65} + 16J_{65}$ holds iff the following 3 equalities hold:
  \begin{equation}
    \label{components_64}
    \begin{cases}
      0^2 + B B^T + 0 = 32 \\
      0 B + BC + B = 16 J_{1,64} \\
      B^T B + C^2 + C = 16 I_{64} + 16J_{64} \\
    \end{cases}
  \end{equation}
  Here $J_{1,64}$ is the matrix of size $1\times64$ consisting of ones.
  
  The first equality of (\ref{components_64}) obviously holds. The second is equivalent to
  \begin{equation}
    \label{partial_sums_C1}
    \sum\limits_{i=1}^{32} C_{ij} =
    \begin{cases}
      15, & j \leqslant 32, \\
      16, & j \geqslant 33. \\
    \end{cases}
  \end{equation}
  Since we also need $AJ_{65} = 32J_{65}$, it follows from (\ref{partial_sums_C1}) that
  \begin{equation}
    \label{partial_sums_C2}
    \sum\limits_{i=33}^{65} C_{ij} = 16\;\;\forall j
  \end{equation}
  Denote $H_{64} = 16J_{64} - B^T B$, the matrix consisting of 4 blocks $15J_{32}$, $16J_{32}$,
  $16J_{32}$, and $16J_{32}$. Also denote $H_{2n} \in M_{2n}(\Z)$ the
  matrix consisting of 4 blocks $15J_{n}$, $16J_{n}$, $16J_{n}$, and $16J_{n}$.

  The last equation of (\ref{components_64}) can be rewritten as
  $$
    C^2 + C = 16I_{64} + H_{64},
  $$
  and now we'll apply the method described in the previous section to this equation. Suppose
  $C$ consists of $4\times4$ circulant blocks of order 16, i.e. it can be `compacted' to $C(x)
  \in M_4(\Z[x]/(x^{16}-1))$ such that
  \begin{equation}
    \label{c_pmod_eq}
    C^2(x) + C(x) \equiv 16I_4 + c_{16}(x)H_4 \pmod{x^{16}-1}.
  \end{equation}
  Let's take this congruence modulo $x-1$ first:
  $$
    C^2(1) + C(1) = 16I_4 + 16H_4
  $$
  As all polynomial coefficients of $C(x)$ must be zeros or ones, elements of $C(1)$ must be
  non-negative integers. By a simple brute-force one can find there are just four such matrices, two
  of them are
  \begin{equation}
    \label{matrices_c1}
    C(1) =
    \begin{pmatrix}
      7 & 8 & 6 & 10 \\
      8 & 7 & 10 & 6 \\
      6 & 10 & 8 & 8 \\
      10 & 6 & 8 & 8 \\
    \end{pmatrix},
    \begin{pmatrix}
      9 & 6 & 7 & 9 \\
      6 & 9 & 9 & 7 \\
      7 & 9 & 6 & 10 \\
      9 & 7 & 10 & 6 \\
    \end{pmatrix},
  \end{equation}
  and the other two are obtained from these two by permuting the third and the fourth rows and
  columns, so we won't consider them. Note that if we successfully expand such matrices into
  $C\in M_{64}(\Z)$, the equalities (\ref{partial_sums_C1}) and (\ref{partial_sums_C2}) will hold.
  
  Next, we consider the congruence (\ref{c_pmod_eq}) modulo $x+1$.
  $$
    C^2(-1) + C(-1) = 16I_4.
  $$
  If $f(x)\in\Z[x]$ is any polynomial with non-negative coefficients and $f(1)=a$ is known, $f(-1)$
  must be a value from $-a$ to $a$ of the same parity, i.e. $f(-1)\in\{-a,\,-a+2,\,-a+4,\,\ldots,\,a-2,\,a\}$.
  So if we know the matrix $C(1)$, we can search for the matrix $C(-1) \equiv C(1) \pmod2$ with
  coefficients less or equal to those of $C(1)$ in absolute values. Using a brute-force approach, we
  find that each matrix in (\ref{matrices_c1}) produces 32 possible matrices $C(-1)$; we can keep only
  10 of them for the first and 10 for the second, the rest will be their permutations. For brevity
  we won't list them all, just one of them obtained from the second matrix (\ref{matrices_c1}):
  \begin{equation}
    \label{example_cm1}
    C(-1) =
    \begin{pmatrix}
      1 & -2 & -3 & -1 \\
      -2 & 1 & -1 & -3 \\
      -3 & -1 & -2 & 2 \\
      -1 & -3 & 2 & -2 \\
    \end{pmatrix}
  \end{equation}

  Next, we consider the congruence (\ref{c_pmod_eq}) modulo $x^2+1$, trying to find possible values of
  $C(i)$.
  $$
    C^2(i) + C(i) = 16I_4.
  $$
  Again, if $f(x)\in\Z[x]$ has non-negative coefficients and we know $f(1)=a$ and $f(-1)=b\equiv
  a\pmod2$, possible values of $f(i)$ have the real part in $\{-\frac{a+b}2,\,-\frac{a+b}2+2,\,
  \ldots,\,\frac{a+b}2-2,\,\frac{a+b}2\}$ and the imaginary part in $\{-\frac{a-b}2,\,-\frac{a-b}2+2
  ,\,\ldots,\,\frac{a-b}2-2,\,\frac{a-b}2\}$. Again, we use brute-force to find possible matrices
  $C(i)$, given $C(1)$ and $C(-1)$. Recall that in order for the final matrix $C$ to be symmetric,
  we need $C^T(x) \equiv C(x^{15}) \pmod{x^{16}-1}$, therefore $C(i) = C^T(-i)$. It turns out there
  are 1422 possible values of $C(i)$ for the first matrix in (\ref{matrices_c1}) and 1224 for the second. One of $C(i)$
  corresponding to (\ref{example_cm1}) is
  \begin{equation}
    \label{example_ci}
    C(i) =
    \begin{pmatrix}
      1 & -2 & -i & -3 i \\
      -2 & 1 & -3 i & -i \\
      i & 3 i & -2 & 2 \\
      3 i & i & 2 & -2 \\
    \end{pmatrix}
  \end{equation}
  Similarly, given $C(1)$, $C(-1)$ and $C(i)$, we find possible matrices $C(\zeta_8)$, where
  $\zeta_8 = \frac{\sqrt2(1+i)}2$, such that
  $$
    C^2(\zeta_8) + C(\zeta_8) = 16I_4
  $$
  and $C^T(\zeta_8) = C(-\zeta_8^3)$. One of possible values corresponding to (\ref{example_cm1}) and
  (\ref{example_ci}) is
  $$
    C(\zeta_8) =
    \begin{pmatrix}
      -1 - 2 \zeta_8 + 2 \zeta_8^3 & 2 \zeta_8^2 & -1 + \zeta_8^2 - \zeta_8^3 & \zeta_8 \\
      -2 \zeta_8^2 & -1 + 2 \zeta_8 - 2 \zeta_8^3 & -\zeta_8 & 1 - \zeta_8^2 - \zeta_8^3 \\
      -1 + \zeta_8 - \zeta_8^2 & \zeta_8^3 & 2 \zeta_8 - 2 \zeta_8^3 & -2 \zeta_8^2 \\
      -\zeta_8^3 & 1 + \zeta_8 + \zeta_8^2 & 2 \zeta_8^2 & -2 \zeta_8 + 2 \zeta_8^3 \\
    \end{pmatrix}
  $$
  Given  $C(1)$, $C(-1)$, $C(i)$ and $C(\zeta_8)$ (the values of $C(x)$ modulo $x-1$, $x+1$, $x^2+1$
  and $x^4+1$), we can calculate $C(x) \pmod{x^8-1}$:
  \begin{multline}
    \label{example_cx}
    \begin{array}{lcrr}
      C_{11}(x) & \equiv & 1 + x^2 + 2 x^3 + 2 x^4 + 2 x^5 + x^6 & \pmod{x^8-1}, \\
      C_{12}(x) & \equiv & x + 2 x^2 + x^3 + x^5 + x^7 & \pmod{x^8-1}, \\
      C_{13}(x) & \equiv & x + x^2 + x^3 + x^4 + x^5 + 2 x^7 & \pmod{x^8-1}, \\
      C_{14}(x) & \equiv & 1 + x + x^2 + 2 x^3 + x^4 + x^6 + 2 x^7 & \pmod{x^8-1}, \\
      C_{22}(x) & \equiv & 1 + 2 x + x^2 + 2 x^4 + x^6 + 2 x^7 & \pmod{x^8-1}, \\
      C_{23}(x) & \equiv & 1 + x^2 + 2 x^3 + x^4 + x^5 + x^6 + 2 x^7 & \pmod{x^8-1}, \\
      C_{24}(x) & \equiv & 1 + x + x^3 + x^5 + x^6 + 2 x^7 & \pmod{x^8-1}, \\
      C_{33}(x) & \equiv & 2 x + x^2 + x^6 + 2 x^7 & \pmod{x^8-1}, \\
      C_{34}(x) & \equiv & 2 + x + x^3 + 2 x^4 + x^5 + 2 x^6 + x^7 & \pmod{x^8-1}, \\
      C_{44}(x) & \equiv & x^2 + 2 x^3 + 2 x^5 + x^6 & \pmod{x^8-1}. \\
    \end{array}
  \end{multline}
  (The rest $C_{ij}(x)$ can be calculated from $C_{ij}(x)\equiv C_{ji}(x^7)\pmod{x^8-1}$). Note
  that this $C(x)$ satisfies the congruence
  $$
    C^2(x) + C(x) \equiv 16I_4 + c_{16}(x) H_4 \pmod{x^8-1}.
  $$
  
  The next step would be to lift the last congruence to $\pmod{x^{16}-1}$, or, equivalently, given
  $C(1)$, $C(-1)$, $C(i)$ and $C(\zeta_8)$ to find possible values of $C(\zeta_{16})$, $\zeta_{16} =
  e^{2\pi i/16}$. Unfortunately, a brute-force doesn't yield any such matrices. Another approach
  could be to introduce a new variable $y$ and try to find a matrix $C(x,\,y)$ such that
  $$
    C^2(x,\,y) + C(x,\,y) \equiv 16I_4 + c_8(x) c_2(y) H_4 \pmod{x^8-1,\,y^2-1}
  $$
  We could then expand $C(x,\,y)$ into $C(x)\in M_8(\Z[x]/(x^8-1))$ and then into $C\in M_{64}(\Z)$.
  However, a brute-force doesn't find any such matrices either.
  
  \section{Last step}
  
  Consider again the calculated matrices $C(x) \pmod {x^8-1}$ such that $ C^2(x) + C(x) \equiv
  16I_4 + 2c_{8}(x) H_4 \pmod{x^8-1}$ (one of them is (\ref{example_cx})). Let's expand them into
  symmetric matrices $D\in M_{32}(\Z)$ consisting of $4\times4$ circulant blocks of order 8
  such that
  \begin{equation}
    \label{d_32}
    D^2 + D = 16I_{32} + 2H_{32}.
  \end{equation}
  
  Suppose now that $C$ consists of $32\times32$ circulant blocks of order 2, therefore it can be
  compacted into $\widetilde{C}(x)\in M_{32}(\Z[x]/(x^2-1))$ such that
  \begin{equation}
    \label{congr_32}
    \widetilde{C}^2(x)+\widetilde{C}(x) \equiv 16I_{32} + c_2(x) H_{32} \pmod{x^2-1}
  \end{equation}
  Furthermore, suppose that $\widetilde{C}(1)$ is one of the matrices $D$ we just found.
  
  Denote $E = \widetilde{C}(0)$, so that $\widetilde{C}(x) = E + (D - E)x$. Then (\ref{congr_32}) is
  equivalent to the two same
  congruences modulo $x-1$ and $x+1$, where the first one holds because of (\ref{d_32}) and the second
  is
  $$
    (2E-D)^2 + (2E-D) = 16I_{32}.
  $$
  After simple transformations we obtain
  \begin{equation}
    \label{de_eq}
    2E^2 - ED - DE + E - D + H_{32} = 0.
  \end{equation}
  We can take this equation modulo 2 to get rid of the square term $E^2$:
  \begin{equation}
    \label{de_eq_mod}
    ED + DE + E + D + H_{32} \equiv 0 \pmod2.
  \end{equation}
  This is a linear equation with respect to $E$. We can solve it as a system of linear equations of $E_{ij}$.
  
  There are several requirements for the final matrix $C$, they impose restrictions on $E$. First,
  $C$ must be symmetric, so must be $E$. Second, all elements of $C$ must be 0 or 1, so must be the
  elements of $E$ and $D-E$ (and therefore, finding $E \pmod2$ is enough to find $E$). Thus,
  $D_{ij}=0$ implies $E_{ij}=0$, $D_{ij}=2$ means $E_{ij}=1$, and for $D_{ij}=1$ $E_{ij}$ can be
  either 0 or 1. Finally, $C$ can have only zeros on the main diagonal because it is an adjacency
  matrix of a graph, thus $E_{ii}=0$. We can iterate over all solutions of (\ref{de_eq_mod}) with such
  restrictions and check if they also satisfy (\ref{de_eq}).
  
  Unfortunately, there are too many solutions, so it may be infeasible to try them all. For example,
  for the matrix $D$ constructed from (\ref{example_cx}) the solutions of (\ref{de_eq_mod}) form an
  affine space of dimension 62, giving $2^{62}\approx4.6\cdot 10^{18}$ possible matrices to check.
  So we need to reduce the set of solutions somehow.
  
  Note that if we multiply some rows and columns of $C(x)$ with the same indices by $x \pmod {x^2-1}$, the
  congruence (\ref{congr_32}) will still hold. In the matrix $C$ that would mean a permutation of some
  pairs of rows and columns, or in the final graph it would be a permutation of some vertices pairs,
  and in the matrix $E$ that corresponds to adding $1\pmod2$ to some rows and columns.
  Therefore, if for every row $i$ we denote $g(i) = \min\{j:\ D_{ij}=1\}$, we may
  assume without loss of generality that $E_{i,g(i)}=0$ if $1\leqslant g(i)<i$. This additional
  restriction can reduce the dimension of solutions of (\ref{de_eq_mod}). For example, for the matrix
  $D$ obtained from (\ref{example_cx}) there would be just $2^{32}\approx4.3\cdot10^9$ solutions, so
  it is now feasible on a modern computer to iterate over all of them, checking whether (\ref{de_eq})
  holds.
  
  And luckily one of such matrices indeed satisfies this equation.
  
  $$
    E =
    \begin{pmatrix}
0 0 0 1 1 1 0 0 0 0 1 0 0 0 0 0 0 0 0 0 0 0 0 1 0 0 0 1 0 0 0 1 \\
0 0 0 0 1 1 1 0 0 0 0 1 0 0 0 0 1 0 0 0 1 1 0 0 1 0 1 0 1 0 0 0 \\
0 0 0 0 0 1 1 1 0 1 0 0 1 0 0 1 0 1 0 0 1 0 1 0 1 1 1 1 0 1 1 0 \\
1 0 0 0 0 0 1 1 1 0 1 0 1 1 1 0 0 0 1 0 1 1 1 1 0 0 1 0 0 1 1 0 \\
1 1 0 0 0 0 1 1 0 0 0 1 0 1 1 0 1 0 0 1 0 1 1 0 1 0 0 1 0 0 1 1 \\
1 1 1 0 0 0 0 1 0 0 0 0 1 0 1 1 0 1 0 0 1 0 1 0 1 1 0 0 1 0 0 1 \\
0 1 1 1 1 0 0 0 1 1 0 1 0 1 0 1 1 1 0 0 0 1 0 1 0 1 0 0 0 1 1 1 \\
0 0 1 1 1 1 0 0 0 1 0 0 1 0 1 0 1 0 0 1 1 0 1 0 0 1 1 1 0 1 1 0 \\
0 0 0 1 0 0 1 0 0 1 1 0 1 0 0 1 1 0 1 1 1 0 1 1 0 1 0 0 0 0 0 1 \\
0 0 1 0 0 0 1 1 1 0 1 1 0 1 0 1 1 1 0 1 1 0 0 0 1 0 0 0 0 0 0 0 \\
1 0 0 1 0 0 0 0 1 1 0 1 1 0 1 0 1 1 0 0 0 1 0 0 1 1 1 0 0 1 0 1 \\
0 1 0 0 1 0 1 0 0 1 1 0 1 1 0 1 0 0 1 1 0 0 1 0 0 0 1 0 1 0 1 0 \\
0 0 1 1 0 1 0 1 1 0 1 1 0 1 1 0 1 1 0 1 0 0 1 1 0 0 1 1 1 0 0 0 \\
0 0 0 1 1 0 1 0 0 1 0 1 1 0 1 0 1 1 1 0 1 1 0 0 0 0 1 1 1 1 0 0 \\
0 0 0 1 1 1 0 1 0 0 1 0 1 1 0 1 0 1 1 0 0 1 0 0 0 1 0 0 0 1 1 0 \\
0 0 1 0 0 1 1 0 1 1 0 1 0 0 1 0 0 0 1 1 0 1 1 1 0 0 1 0 0 0 1 1 \\
0 1 0 0 1 0 1 1 1 1 1 0 1 1 0 0 0 1 0 0 0 0 0 1 1 1 0 1 1 0 1 0 \\
0 0 1 0 0 1 1 0 0 1 1 0 1 1 1 0 1 0 1 0 0 0 0 1 1 1 0 0 1 1 0 1 \\
0 0 0 1 0 0 0 0 1 0 0 1 0 1 1 1 0 1 0 1 0 0 0 0 1 1 1 1 0 0 1 0 \\
0 0 0 0 1 0 0 1 1 1 0 1 1 0 0 1 0 0 1 0 1 1 0 0 1 1 0 1 1 0 0 1 \\
0 1 1 1 0 1 0 1 1 1 0 0 0 1 0 0 0 0 0 1 0 1 1 0 1 1 1 0 1 1 0 0 \\
0 1 0 1 1 0 1 0 0 0 1 0 0 1 1 1 0 0 0 1 1 0 1 1 1 1 1 1 0 1 0 0 \\
0 0 1 1 1 1 0 1 1 0 0 1 1 0 0 1 0 0 0 0 1 1 0 1 0 0 1 1 1 0 1 1 \\
1 0 0 1 0 0 1 0 1 0 0 0 1 0 0 1 1 1 0 0 0 1 1 0 0 0 0 1 1 1 1 1 \\
0 1 1 0 1 1 0 0 0 1 1 0 0 0 0 0 1 1 1 1 1 1 0 0 0 0 0 1 0 1 1 0 \\
0 0 1 0 0 1 1 1 1 0 1 0 0 0 1 0 1 1 1 1 1 1 0 0 0 0 0 0 1 0 1 1 \\
0 1 1 1 0 0 0 1 0 0 1 1 1 1 0 1 0 0 1 0 1 1 1 0 0 0 0 0 1 1 0 1 \\
1 0 1 0 1 0 0 1 0 0 0 0 1 1 0 0 1 0 1 1 0 1 1 1 1 0 0 0 0 1 1 0 \\
0 1 0 0 0 1 0 0 0 0 0 1 1 1 0 0 1 1 0 1 1 0 1 1 0 1 1 0 0 0 1 1 \\
0 0 1 1 0 0 1 1 0 0 1 0 0 1 1 0 0 1 0 0 1 1 0 1 1 0 1 1 0 0 0 1 \\
0 0 1 1 1 0 1 1 0 0 0 1 0 0 1 1 1 0 1 0 0 0 1 1 1 1 0 1 1 0 0 0 \\
1 0 0 0 1 1 1 0 1 0 1 0 0 0 0 1 0 1 0 1 0 0 1 1 0 1 1 0 1 1 0 0 \\
    \end{pmatrix}
  $$
  
  From this matrix and from $D$ we can construct $C(x)\in M_{32}(\Z[x]/(x^2-1))$, then we can expand $C(x)$
  to $C\in M_{64}(\Z)$. Then by adding a row and a column of $B$ (32 ones followed by
  32 zeros) we finally obtain the adjacent matrix $A$ of a SRG with the parameters $(65,\,32,\,15,\,16)$,
  we also write this matrix below. One can easily check that it consists only of zeros and ones, that it
  is symmetric, that each row has 32 ones, and each two rows $i$, $j$ have 15 or 16 common ones,
  depending on $A_{ij}$.

  Note that we can also construct a symmetric conference matrix of order 66 from it by replacing
  all $1$ with $-1$, then all $0$ with 1, keeping zeros on the main diagonal and adding a row and a
  column of ones. A square of such matrix is equal to $65I_{66}$.
  
  One can also compute the automorphism group of this SRG. It consists of 32 elements, thus the
  graph is not vertex-transitive. In fact, its vertex orbits are $\{1\}$, $\{2,\,3,\,\ldots,\,33\}$,
  and $\{34,\,35,\,\ldots,\,65\}$.
  
  Here is the adjacency matrix of the SRG:

\phantom{aaa}\\
$01111111111111111111111111111111100000000000000000000000000000000$ \\
$10100011111110100000111010001000100010101010100110101011101000111$ \\
$11000101111111000001011100010001000101010101000111010101110001011$ \\
$10001000111111101010001110100010011000101101001001101100111010001$ \\
$10010001011111110100010111000100011001010010110001110011011100010$ \\
$10100010001111111001000011101001000110001100110011011101001111000$ \\
$11000100010111111000100101110000100110010011001100111010110110100$ \\
$11101000100011111100010001011100001001100101010100001110101101101$ \\
$11110001000101111010001000111010010001100010101010010111010011110$ \\
$11111010001001011000100100010110110010011001010011000011101011011$ \\
$11111100010000111001000010001111001100011000101100100101110100111$ \\
$11111110100010010010001001000101101100100110010011110000111010110$ \\
$11111111000100001100010000100011110011000110001101101001011101001$ \\
$10111111110000100111000100010001010100101001100100111010001111010$ \\
$11011111101001000110100010001000101011010001100011011100010110101$ \\
$10001111111100001011101001000100010010110100011000110111000101101$ \\
$10010111111010010101110000100010001101001010011001001110100011110$ \\
$10001001000011101011110001100011110001011100110110110000100010111$ \\
$10010000100101110101101001100101101000111011001111001001000101011$ \\
$10100100001001011110111100011001011100010110101011101010001000101$ \\
$11000010010000111111011010011000111010001111010101110100010001010$ \\
$11101001000010001101101111000110010110100011101011011100100100010$ \\
$11110000100100010011110110100110001111000101110100111011000010001$ \\
$10111010010001000001011011110001101011110000111010101110110001000$ \\
$11011100001000100000111101101001110101101001011101010111001000100$ \\
$10001111000100010110010110111100010100111010010110001101110010001$ \\
$10010110100010001110001111011010001011011100001110010011101100010$ \\
$10100011110001000001100101101110111101001111000010100101011100100$ \\
$11000101101000100001100011110111011010110110100101000010111011000$ \\
$10001001011100010010011001011011101111001011101000010000101111001$ \\
$10010000111010001100011000111101110110110101110000001001010110110$ \\
$10100100001111000111000110001110100011110011011100100100001011110$ \\
$11000010010110100110100110010111000101101100111011000010010101101$ \\
$00011000110011010101110011011010000110100000001111110001011011101$ \\
$00011001001100101011101100111100000111000000010111101000111101110$ \\
$00100110001101001001011011010110111001101000000101011010010110111$ \\
$01000110010010110000111100101111011001110000000010111100001111011$ \\
$00101001100010101100001110110101101110011010000001110111000011101$ \\
$01010001100101010010010111001011110110011100000001101110100101110$ \\
$00101010011000110111000101101011000011100111000001011011110000111$ \\
$01010100011001001110100011110100100101100110100000111101101001011$ \\
$00110101000110010101101000111010100000111001110001110110111100001$ \\
$01001010100110001011110001011101000001011001101001101111011010010$ \\
$00110011010001100010111010010111000000010110011101011101101110100$ \\
$01001100101001100101011100001110100000001110011010111011110111000$ \\
$00001101010100011100101111000011101000000101100110001111011011110$ \\
$00010010101010011011010110100101110000000011100110010110111101101$ \\
$01100011001011000110101011101001011100000001011000100011110111011$ \\
$01100100110100100111010101110000111010000000111001000101101110111$ \\
$00111100010110101011110010001000111101110111000010000011100111000$ \\
$01011010001111010101101100010001011011101110100100000101100110100$ \\
$00101110100101110100111010100100010111011101101000000000111001110$ \\
$01010111000011101011011101000010001110111011110000000001011001101$ \\
$00110101101000111000110111010001000011101111011010100000010110011$ \\
$01001011110001011001001110101000100101110110111101000000001110011$ \\
$01101100111010010010001011110010010001011011110111101000000101100$ \\
$01110011011100001100010101101100001000111101101111110000000011100$ \\
$00111010101110100000100101011010111100010110111100011100000001011$ \\
$01011101010111000001000010111101011010001111011010011010000000111$ \\
$00001111001011110010010000110110101110100101101111100111000000010$ \\
$00010110110101101100001001001111010111000011110111100110100000001$ \\
$00100101110011011010100100001101111011101000111101011001110000000$ \\
$01000011101100111101000010010011111101110001011010111001101000000$ \\
$01101000111101001110110000100011001110111010010110010110011100000$ \\
$01110001011010110111001001000100110111011100001110001110011010000$ \\

\end{document}